\documentstyle[12pt]{article}
\normalbaselines
\begin{document}

\vbox {\vspace{6mm}}

\begin{center}
{\large \bf Tensor product in symmetric function spaces}\\[2mm]
S.V.   Astashkin\\
{\it Samara State University,443011 Samara 11, Russia}\\[5mm]
\end{center}

\begin{abstract}
A concept of multiplicator of symmetric function space concerning to
projective tensor product is introduced and studied. This allows us to
obtain some concrete results. In particular, the well-known theorem of
R.O'Neil about the boundedness of tensor product in the Lorentz spaces $L_{pq}$ is discussed. \vspace{\baselineskip}
\end{abstract}
\begin{center}
0. {\em Introduction}
\end{center}

Let $x=x(s)$ and $y=y(t)$ are measurable functions on
$I=[0,1]$. We define following bilinear operator:
$$B(x,y)(s,t)=(x\otimes{y})(s,t)=x(s)y(t),(s,t)\in{I\times{I}}$$
If $X,Y,Z$ are symmetric function spaces, then the boundedness of $B$
from $X\times{Y}$ into $Z(I\times{I})$ is equivalent to the continuity of
the embedding
$$
{X\otimes{Y}}\subset{Z}$$
where $X\otimes{Y}$ denotes the projective tensor product of spaces $X$
and $Y$ [ 1, p.51 ].

The importance of study of operator $B$ follows mainly from two facts : this
operator has various applications in the theory of symmetric spaces
( see, for example, [2,p.169-171], [3], [4] ) and the problem of
boundedness of tensor product is equivalent to the stability problem of
integral operators ( [ 1, p. 50 ], [ 5 ] ).

The problem of boundedness of the operator  $B$ in certain families
of Banach spaces  ( for example, Lorentz, Marcinkiewicz and Orlicz spaces )
has attracted the attention of number of authors ( see [ 1 ], [ 6 - 11 ] ).
In particular, R.O'Neil proved following theorem about the
boundedness of tensor product in the spaces $L_{pq}$ [ 9 ].

Let us remind that the space $L_{pq}(1<p<\infty, 1\le{q}\le{\infty})$
consists of all measurable on $I$ functions $x=x(t)$ for which
$$
|| x ||_{pq} = \left\{
\begin{array}{ll}
\left\{\int\limits_0^1\left(x^*(t)t^{1/p}\right)^q\,{dt}/t\right\}^{1/q}, & \mbox{if } 1\le{q}<\infty\\
\sup\limits_{0<t\le{1}}\left(x^*(t)t^{1/p}\right),                        & \mbox{if } q=\infty
\end{array} \right. $$
is finite. Although the functional $\|x\|_{pq}$ is not subadditive, however it is equivalent to
the norm $||x||_{pq}'\,=\,||x^{**}||_{pq}$. where $x^{**}(t)={1/t}{\int_0^t{x^*(s)}\,ds}$
and by $x^*(s)$ is denoted the left continuous non-increasing rearrangement
of $|x(s)|$.

Let us notice that ${L_{pr}}\subset{L_{pq}} ( r\le q ),{}L_{pp}=L_p$ and
usual norm of $L_p$ $||x||_p=||x||_{pp}$. These spaces arise naturally in
interpolation theory ( see, for example, [ 12 ] ).

{\bf Theorem ( O'Neil )}. Let $1<p<\infty, 1\le{q,r,s}\le\infty$. The boundedness
of the operator $B$ from $L_{pr}\times{L_{pq}}$ into $L_{ps}(I\times{I})$ is
equivalent to the conditions :
$$
1)\max(q,r)\le{s}$$
and
$$
2){\frac1p}+{\frac1s}\le{\frac1q+\frac1r}$$.

We consider following problems connected with the exactness of this theorem :

a ) Let us fix $1<p<\infty,1\le{r}\le{s}\le{\infty}$ and $q^*$ is the
maximum of all $q\ge{1}$ for which the conditions of theorem  are satisfied
( if those exist ). Is the space $L_{p{q}^*}$ the largest among symmetric spaces
$E$ such that
$$
B: { L_{pr}\times{E}}\rightarrow {L_{ps}(I\times{I})}\,\,?$$

b ) Let us $1<p<\infty,1\le{r}\le{q}\le{\infty}$ and $s^*$ is the minimum
of all $s\ge{1}$ for which the conditions of theorem  are satisfied ( if
those exist ). Is the space $L_{p{s^*}}$ the smallest among symmetric spaces
$E$ such that
$$
B:{L_{pr}\times{L_{pq}}}\rightarrow{E(I\times{I})}\,\,?$$

The answers on these questions are positive, of course, in case
$r\le p$, since B bounded from $L_{pr}\times{L_{pr}}$ into $L_{pr}(I\times{I})\,(r\le{p})$
and from $L_{p}\times{L_{pq}}$ into $L_{pq}\,(p\le{q})$ ([9],[11]).

In this paper we shall study tensor product in general symmetric spaces.
It will allow us to solve in part the problems a) and b).
$$\,\,$$

The paper is organized in the following manner. At first, we introduce
a concept of the multiplicator $M(E)$ of a symmetric space $E$ concerning to
tensor product. The relations defining the fundamental function and the
norm of dilation operator in $M(E)$ will be obtained. Moreover, the upper
and lower estimates of $M(E)$ allow to find it in some cases. In particular,
we shall get the positive answer on the question a) in case $p<r=s$ ( see
the corollary 1.12 ). Then we shall consider a quantity connected with the
boundedness of tensor product in symmetric spaces.

Finally we study the boundedness of operator $B$ in the spaces $L_{pq}$
using the family of Lorentz-Zygmund spaces $L_{pq}(\log\,L)^\alpha$ that
involve, in particular, spaces $L_{pq}$. The main result --- theorem 2.1 ---
is similar to O'Neil theorem. First of all, it allows to define
more exactly values of $B$ on $L_{pr}\times{L_{pq}}$ for $p<r\le{q},\,q^{-1}+r^{-1}-p^{-1}\ge{0}.$
Namely, we shall find ( see corollary 2.7 ) the symmetric space
${E=E(p,r,q)}\subset{L_{ps}},E\ne{L_{ps}},\,s^{-1}=q^{-1}+r^{-1}-p^{-1}$ such that
$$
B: L_{pr}\times{L_{pq}}\rightarrow{E(p,r,q)\,(I\times{I})}$$

The last means that the answer on the problem b) is negative in this case.
Besides, O'Neil theorem does not describe the set $B\,(\,L_{pr},\,L_{pq}\,)$,
if $r^{-1}+q^{-1}-p^{-1}<0$ ( a family $L_{pq}$ is too narrow for study
of this problem ). Using the spaces $L_{pq}(\log\,L)^\alpha$ we solve this problem
to within any positive power of logarithm
( see theorems 2.1 and 2.4 ). The complex method of interpolation is used
here ( the definitions of interpolation theory see in [12] or [13] ).\newpage
\begin{center}
1. {\em Thr multiplicator of symmetric space concerning to tensor product}
\end{center}

If $z=z(w)$ is a measurable function on $I$ or $I\times{I},\,\mu$ is the usual Lebesgue measure,
then the distribution function of $|z(s)|$ is $n_{z}\,(\tau)=\mu\{w:\,|z(w)|>\tau\,\}\,\,(\,\tau>0\,).$
The function $x(t)$ is equimeasurable with $y(t)$ if $n_x(\tau)=n_y(\tau)\,(\tau>0).$

 Let us remind that a Banach space $E$ of measurable functions defined
on $I$ is said to be symmetric if these conditions are satisfied: a) $y\in{E}$
and $|x(t)|\le{|y(t)|}$ imply that $x\in{E}$ and $||x||\le{||y||};$ b) if
$y\in{E}$ and the function $x(t)$ is equimeasurable with $y(t)$, then
$x\in{E}$ and $||x||=||y||.$

Let $E$ be a
symmetric space on $I=\,[0,1]$. We shall denote by $M(E)$ the set of
all measurable functions $x=x(s)\,(\,s\in{I}\,)$ for which  $x\otimes{y}\in{E(I\times{I})}$
with arbitrary $y\in{E}$. Then $M(E)$ is symmetric space on $I$ concerning to the norm
$$
||x||\,=\, sup\,\{\,||x\otimes{y}||_{E(I\times{I})}:\,||y||_E\,\le1\,\} $$
Obvios, $M(E)\,\subset{E}$.

{\bf Example 1.1}. Let us find $M(E)$, if $E$ is Lorentz space $\Lambda(\,\phi\,)$,
where positive concave function $\phi(u)$ increases on $(0,1]$.
This space consists of all measurable functions $x=x(s)$
for which $$ ||x||_{\Lambda(\,\phi\,)}\,=\int_0^1 x^*(s)\,d\phi(s)\,\,<{\infty}$$
In particular, $L_{p1}\,=\,\Lambda(\,t^{1/p}\,\,)\,\,(\,1<p<\infty\,)$.

For arbitrary $e\subset{I},\,x\in{\Lambda(\,\phi\,)}$ functions
$\chi_e\otimes{x}$ and $\sigma_{\mu(e)}x$ are
equimeasurable, where $\chi_e\,(s)=1\,(s\in{e}),\,\chi_e\,(s)=0\,(s\notin{e})$ and
the dilation operator $\sigma_{t}y(u)=y(\,u/t\,)\chi_{[0,1]}\,(\,u/t\,)\,(t>0).$
The rearrangements of equimeasurable functions are equal a.e.. Therefore,
$$
||\chi_e\otimes{x}||_{\Lambda(\,\phi\,)}\,=\int_0^1(\sigma_{\mu(e)}\,x)^*(u)\,d\phi(u)=\,\int_0^{\mu(e)}x^*(u/\mu(e))\,d\phi(u)=\,\int_0^1x^*(v)\,d\phi(\mu(e)v)$$
Since $||\chi_e||_{\Lambda(\,\phi\,)}=\phi(\mu(e))$, then we have in view of [ 12, p.151 ]
$$
||x||_{M(\,\Lambda(\,\phi\,)\,)}\,\le\,2\sup_{0<t\le1}\int_0^1x^*(v)\,d{\left [\frac{\phi(tv)}{\phi(t)}\right]}.$$

Introduce the notation:
$$
{\cal{M}}_\phi(v)\,=\,\sup\,\left\{\,\frac{\phi(tv)}{\phi(t)},\,\,\,0<t\le\min(1,1/v)\right\}.$$
Using last inequality, we obtain
$$
\Lambda(\,{\cal{M}}_\phi\,)\subset{M(\,\Lambda(\,\phi\,)\,)}$$

Assume in addition that
$$
{\cal{M}}_{\phi}(v)\,=\,\lim_{t\to t_0}\frac{\phi(tv)}{\phi(t)},\eqno{(1)}$$
where $t_0\in{[0,1]}$ does not depend on $v\in{[0,1]}$. In this case
$$
||x||_{\Lambda(\,{\cal{M}}_\phi\,)}\,=\,\lim_{t\to {t_o}}\,\int_0^1x^*(v)\,d{\left [\frac{\phi(tv)}{\phi(t)}\right]}\,\,=\,\,\lim_{t\to{ t_0}}\frac{||\chi_{(0,t)}\otimes{x}||_{\Lambda(\,\phi\,)}}{\phi(t)}\,\,\le\,||x||_{M(\,\Lambda(\,\phi\,)\,)}.$$
 
Consequently, if (1) holds, then
$$
M(\,\Lambda(\,\phi\,)\,)\,=\,\Lambda(\,{\cal{M}}_\phi\,).$$\vspace{\baselineskip}

Let again $E$ be arbitrary symmetric space on $I$. We shall obtain the relations
connecting the fundamental functions and the norms of dilation operators
of spaces $E$ and $M(E).$
\def\ab{\phi_{M(E)}(t)}
\newcommand{\bc}{||\sigma_t||_{{E}\to{E}}}
\def\cd{||\sigma_t||_{M(E)\to M(E)}}
\def\mn{||\sigma_{1/t}\,||_{E\to{E}}}
\def\ef{||\sigma_{1/t}||_{M(E)\to{M(E)}}}
\def\fg{x\otimes{y}}

In the first place,
$$
\ab\,=\,\bc,\,\,\,0<t\le1,\eqno{(2)}$$
where $\phi_X(t)$ denotes the fundamental function of symmetric space $X$:
$\phi_X(t)=\,||\chi_{(0,t)}||_X$. Indeed, since the functions $\chi_{(0,t)}\otimes{y}$
and $\sigma_ty\,\,(\,t\in{[0,1]}\,)$ are equimeasurable, then
$$
\ab\,=\,\sup_{||y||\le1}||\chi_{(0,t)}\otimes{y}||_E\,\,=\,\,\sup_{||y||\le1}||\sigma_ty||_E\,\,=\,\,\bc.$$

{\bf Theorem 1.2}. For arbitrary symmetric space $E$ on $I$ :
$$
\cd\,\,=\,\,\bc,\,\,\,0<t\le1\eqno{(3)}$$
$$
{\mn}^{-1}\,\,\le\,\,\cd\,\,\le\,\,\bc\,,\,\,\,t>1\eqno{(4)}$$
{\bf Proof}. First of all, if $x(u)$ and $y(v)$ are measurable functions on $I$ and
$$
supp\,x\,=\,\{u:\,x(u)\ne0\}\,\subset\,{[0,\min(1,t^{-1})]},$$
then for $t,s>0$
 $$n_{\sigma_t\fg}(s)=\,\int_0^1\,\mu\left\{\,u:\,|\sigma_tx(u)|\,>\,\frac s{|y(v)|}\,\right\}\,dv=t\int_0^1\,\mu\left\{\,u:\,|x(u)|\,>\frac s{|y(v)|}\,\right\}\,dv=t\,\,\,n_{\fg}(s).$$

Consequently,
$$
 (\,\sigma_{t}{\fg}\,)^*\,(s)\,=\,\sigma_{t}(\fg)^*\,(s)\,\,\,\,0\le{s}\le{1},$$
and thus
$$
||\sigma_t\fg||_E\,=\,||\sigma_t(\fg)^*||_E\,\le\,\bc\,||\fg||_E\eqno{(5)}$$
If $0<t\le1$, then we get from (5) by definition of norm in $M(E)$
$$
\cd\,\,\le\,\bc\eqno{(6)}$$
In case $t>1$ we define the function $\tilde{x}(u)=x(u)\chi_{(0,1/t]}(u).$
For it (5) is satisfied and, besides, the functions $\sigma_t{\tilde{x}}$ and
 $\sigma_t{x}$ are equimeasurable. Therefore,
$$
||\sigma_t\fg||_E\,\,=\,||\sigma_t{\tilde{x}}\otimes{y}||_E\,\,\le\,\bc\,||\tilde{x}\otimes{y}||_E\,\,\le\,\bc\,||\fg||_E,$$
and thus (6) is satisfied for $t>1$.

If $0<t\le1,$ then in view (2) and by the quality $||\chi_{(0,1)}||_{M(E)}\,=\,1$
$$
\bc\,\,=\,\ab\,=||\sigma_t{\chi_{(0,1)}}||_{M(E)}\,\,\le\,\cd$$

Using the submultiplicativity of function $f(t)=\,\cd\,\,(t>0)$ and the inequality
(6) we get in case $t>1$:
$$
\cd\,\,\ge\,\,{\ef}^{-1}\,\,\ge\,\,{\mn}^{-1}$$

The theorem is proved.
$$\,\,$$

The Boyd indices of symmetric space $E$ are defined as
$$
\alpha_{E}\,=\,\lim_{t\to0}\frac{\bc}{\ln{t}}\,\,,\,\,\beta_{E}\,=\,\lim_{t\to{\infty}}\frac{\bc}{\ln{t}}\,.$$
From theorem 1.2 follows

{\bf Corollary 1.3}. If $E$ is any symmetric space on $I$, then
$$
0\,\le\alpha_E\,=\,\alpha_{M(E)}\,\le\,\beta_{M(E)}\,\le\,\beta_E\,\le\,1.$$

{\bf Remark 1.4}. Left-hand side of the inequality (4) cannot be improved,
in general. We shall show that there exists a symmetric space $E$ on $I$ such that
$$
\cd\,\mn\,\,=\,1\,(t\ge1)\eqno{(7)}$$
and
$$
\lim_{t\to{\infty}}\frac{\cd}{\bc}\,\,=\,0\eqno{(8)}$$

Let $\phi_{\alpha}(s)=\,s^{\alpha}\ln^{-1}(C/s),\,E=\,\Lambda(\,\phi_{\alpha}\,)\,(0<\alpha<1).$
If $C>\exp\{1/(1-\alpha)\},$ then the function $\phi_{\alpha}$
is concave on $(0,1]$. The function
$f_t(s)=\,\ln(C/s)\,\ln^{-1}[C/(st)]\,\,(0<s\le\min(1,1/t)\,)$
decreases, if $0<t\le1$, and increases, if $t>1.$ Therefore, for $0<t\le1$
\def \gh {\phi_{\alpha}}
$$
{\cal{M}}_{\gh}(t)=\, t^{\alpha}\lim_{s\to 0+}f_t(s)=\,t^{\alpha}.$$
Hence, in particular, it follows that the condition (1) is fulfilled for $\gh(s)$
and consequently as it was showed in example 1.1
$$
M(\,\Lambda(\,\gh\,)\,)\,=\,\Lambda(\,{\cal{M}}_{\gh}\,)\,=\,\Lambda(\,t^{\alpha}\,)$$

In the same time, for $t>1$
$$
{\cal{M}}_{\gh}(t)\,=\,t^{\alpha}f_{t}(1/t)\,=\,t^{\alpha}{\ln(Ct)}\,{\ln^{-1}{C}}$$
Since $||\sigma_t||_{\Lambda(\phi)\to {\Lambda(\phi)}}\,\,=\,\,{\cal{M}}_{\phi}(t)\,\,\,$
[ 12, p. 134 ], then in our case
$$
\bc\,\,=\,\,{\cal{M}}_{\gh}\,(t),\,\,\,\cd\,\,=\,\,t^{\alpha}\,\,\,(t>0)$$
As a result we obtain (7) and (8).\vspace{\baselineskip}

Continuing the study of multiplicator of a symmetric space, we shall find its upper
and lower estimates in general case.

{\bf Theorem 1.5(upper estimate)}. Suppose that $E$ is a symmetric space on $I$,
$p=\,1/\alpha_E,$ where $\alpha_E$ is lower Boyd index of space $E$. Then
$$
M(E)\,\subset\,{L_p}$$
and the constant of this embedding does not depend from a space $E$.
\def \hi {\sum\nolimits_{i=1}^m\alpha_{i}\chi_{(\frac{i-1}{m},\frac{i}{m}]}}
\def \ik {\alpha}
\def \km {\left(\sum_{k=1}^m{\alpha_k}^p\right)^{1/p}}
\def \op {\epsilon}
\def \pq {||\sigma_t||}
\def \qr {\sum_{k=1}^m{\alpha_{k}y_k}}
\def \rs {{\cal{M}}_{\phi_E}}

{\bf Proof}. We consider the nonnegative simple functions of type:
$$
x(s)=\,\,\hi(s),\,\,\ik_{i}\ge 0,\,m\in{\bf{N}}\eqno{(9)}$$
For any $y\in{E}$ function $\fg$ is equimeasurable with the function
$$
z_{x,y}(t)\,=\,\,\qr(t),$$
where $y_{k}\in{E}$ are disjointly supported functions having all the same
distribution function and
$$
n_{y_k}\,(s)\,=\,\frac{1}{m}n_{y}(s),\,\,s>0,\,k=1,2,..,m\eqno{(10)}$$
We shall show that for every $\op>0$ there exists $y\in{E}$ ( or
$y_k,\,\,k=1,2,..,m$ ) for which the equalities (10) hold and
$$
||\fg||_{E}\,=\,||z_{x,y}||_{E}\,\,\ge{\frac{1-\op}{1+\op}\,||x||_{p}\,||y||_{E}},\eqno{(11)}$$
where $x$ is arbitrary function of type (9) and $p=1/\ik_{E}.$

Indeed, in view of [2, p. 141] ( see also [14] ) one can choose disjointly
supported equimeasurable functions $z_k\in{E}\,(k=1,2,..,m)$ such that for
$\ik_{k}\ge{0}$
$$
(1-\op)\km\,\le{||\sum_{k=1}^m{\ik_{k}z_{k}}||_{E}}\,\le{(1+\op)\km}\eqno{(12)}$$
Denote
$$
y_{k}=\frac{z_{k}}{m}\,(k=1,2,..,m),\,\,y=\,\frac{1}{m}\,\sum_{k=1}^m{z_{k}}.$$
These functions satisfy all necessary assumptions, in particular, the
equalities (10). We shall prove (11).

From (12) ( if we take $\ik_{k}=1$ for $k=1,2,..,m$ ) it follows that
$$
||y||_{E}\,\le{(1+\op)m^{1/p-1}},$$
and therefore
\begin {eqnarray*}
||\qr||_{E}\,=\,\frac{1}{m}||\sum_{k=1}^m{\ik_{k}z_{k}}||_{E}\,\ge{\frac{1-\op}{m}\,\km}\ge\\
     \ge{\frac{1-\op}{1+\op}m^{-1/p}\,||y||_{E}\,\km}\,=\frac{1-\op}{1+\op}\,||x||_{p}\,||y||_{E}.
\end {eqnarray*}

Since $\op>0$ is arbitrary, then the following inequality can be obtained
for functions $x(s)$ of type (9):
$$
||x||_{M(E)}\,\ge{||x||_{p}}\eqno{(13)}$$

Let us extend last inequality on countably-value dyadic functions
$$
x(s)\,=\,\sum_{k=1}^{\infty}{\ik_{k}\,\chi_{(2^{-k},2^{-k+1}]}\,(s)},\,\,\ik_{k}\ge0.$$
Functions $x_{n}(s)=x(s)\chi_{(2^{-n},1]}(s)$ belong to the class (9),
 $x_{n}(s)\le{x(s)}\,(n=1,2,..).$ Consequently,
$$
||x||_{M(E)}\,\ge{||x_{n}||_{M(E)}}\,\ge{||x_{n}||_{p}}\,\,(n=1,2,..).$$
Passing to the limit as $n\to{\infty}$ we obtain (13).

If now $y\in{E}$ is arbitrary, then we define the function
$$
x_{y}(s)\,=\,\sum_{k=1}^{\infty}{y^*(\,2^{-k+1}\,)\,\chi_{(2^{-k},2^{-k+1}]}\,(s)}.$$
It is clear that $x_{y}\le{y^*}\le{\sigma_{2}(x_{y})}.$ Since for any symmetric space $X$
$||\sigma_{t}||_{X\to{X}}\,\le{\max(1,t)}$ [12, p. 133], then
$$
||y||_{M(E)}\,=\,||y^*||_{M(E)}\,\ge{||x_{y}||_{M(E)}}\,\ge{||x_{y}||_{p}}\,\ge\,\frac{1}{2}||y||_{p}.$$

The theorem is proved.
$$\,\,$$

{\bf Theorem 1.6 (lower estimate)}. For arbitrary symmetric space $E$ on $I$
$$
M(E)\,\,\supset\,\,{\Lambda(\,\psi\,)},$$
where $\psi(t)=\bc\,\,(\,0<t\le1\,).$

Moreover, the constant of embedding does not depend from a space $E.$

{\bf Proof}. We fix $y\in{E}$ and consider the operator $T_{y}x=\fg.$ If
$x=\chi_{e},\,e\subset{[0,1]},\,\mu(e)=t,$then
$$
||T_{y}\chi_{e}||_{E}\,=\,||\sigma_ty||_E\le{\bc\,||y||_E}\,=\,\psi(t)||y||_E.$$
In view of [12, p. 151] the following inequality is true for every $x\in{\Lambda(\,\psi\,)}:$
$$
||\fg||_E\,=\,||T_yx||_E\,\le\,{2||y||_E\,||x||_{\Lambda(\psi)}}.$$

The theorem is proved.\vspace{\baselineskip}

Consider some corollaries of obtained theorems.

{\bf Corollary 1.7}. The multiplicator $M(E)=L_{\infty}$ iff the lower Boyd index
of symmetric space $E$  $\ik_E=0$.

{\bf Proof}. If $\ik_E=0$, then $M(E)=L_{\infty}$ in view of theorem 1.5.
If $\ik_E>0$, then by theorem 1.6 $M(E)\supset{\Lambda(\psi)}$ and
$\psi(t)=||\sigma_t||_{E\to{E}}\,\to{0}\,(t\to{0+})$. Therefore $\Lambda(\psi)\neq{L_{\infty}}.$
Since $\Lambda(\psi)\supset{L_{\infty}}${\newline} [ 12,p.124 ] then, all the more,
$M(E)\neq{L_{\infty}}.$

{\bf Corollary 1.8}. If $B:\,\,E\times{E}\to{E(I\times{I})}$, then $E\subset{L_{1/{\ik_E}}}$.\vspace{\baselineskip}
$$\,\,$$

If $E$ is a Banach space of measurable functions defined on $I$, then dual space $E'$ consists of all measurable on $I$
functions $y=y(t)$ for which
$$
||y||_{E'}\,=\,\sup\{\int_0^1x(t)y(t)dt:\,||x||_{E}\le1\}\,<\,\infty.$$

The norm of symmetric space $E$ is called order subcontinuous
if from monotone convergence $x_n\uparrow{x}$ a.e.
$(x_n,x\in{E})$ it follows: $||x_n||_{E}\to{||x||_{E}}.$

{\bf Corollary 1.9}. Let $E$ be a symmetric space with order subcontinuous norm. If
$B:\,\,E\times{E}\to{E(I\times{I})}$ and $B:\,\,E'\times{E'}\to{E'(I\times{I})}$,
then there exists $p\in{[1,\infty]}$
such that $E=L_p.$

{\bf Proof}. By corollary 1.8 $E\subset{L_{1/\ik_E}}.$ For dual spaces the inverse
embedding is satisfied:
$$
E'\,\supset{(L_{1/\ik_E})'}\,=\,L_{1/{1-\ik_E}}\,\supset{L_{1/\beta_{E'}}}\,\supset{L_{1/\ik_{E'}}},$$
since $1-\ik_E\ge{\beta_{E'}}\ge{\ik_{E'}}$ [12, p. 144].

On the other hand, now in view of corollary 1.8 $E'\subset{L_{1/\ik_{E'}}}$
and thus $E'=L_{1/\ik_{E'}}.$

As $L_{\infty}\subset{E}$ and $E$ is the subspace of
$E''$ [15, p.255] we obtain, finally, for $p=1/(1-\ik_{E'})$ the equality:
$$
E\,\, =\,\,E''\,\,=\,\,L_p.$$\vspace{\baselineskip}

If $\phi_E(s)$ is the fundamental function of a symmetric space $E$, then
${\rs}(t)\le{\bc}$. A space $E$ is called space of fundamental type provided     									
$$
\bc\,\,\le\,\,{C\rs\,(t)},$$
where $C$ does not depend on $t>0$. Let us notice that all most important in
applications symmetric spaces are spaces of fundamental type.

{\bf Corollary 1.10}. If $E$ is a symmetric space of fundamental type, then $M(E)\supset{\Lambda(\,\rs\,)}$\vspace{\baselineskip}

Let us remind that a Banach space $E$ is interpolation space between the spaces
$E_0$ and $E_1$ if $E_0\cap{E_1}\subset{E}\subset{E_0+E_1}$ and from
boundedness of arbitrary linear operator in $E_0$ and $E_1$ follows its
boundedness in $E$. In this case there exists $C>0$ such that
$$
||T||_{E\to{E}}\,\le{C\max_{i=0,1} {||T||_{E_i\to{E_i}}}}$$
where $C$ does not depend on operator $T$[12].

{\bf Theorem 1.11}. Let symmetric space $E$ be an interpolation space between the spaces
$L_p$ and $L_{p\infty}$ by some $p\in{(1,\infty)}.$ Then $M(E)=L_p.$

{\bf Proof}. Since $L_p\subset{E}\subset{L_{p\infty}},$ then, obvious, there
exist $C_1>0$ and $C_2>0$ for which
$$
C_1t^{1/p}\le{\phi_E(t)}\le{C_2t^{1/p}}\,\,(0\le{t}\le{1}).$$
Hence, in particular,
$$
\bc\,\,\ge\,\,C_3t^{1/p}\,\,(t>0).$$

On the other hand, by condition of theorem
$$
\bc\,\,\le{C_4\max\left(||{\sigma}_t||_{L_p\to{L_p}}\,,\,||{\sigma}_t||_{L_{p\infty}\to{L_{p\infty}}}\right)}\,=\,C_4t^{1/p}.$$
Using obtained inequalities, we get in view of definition of Boyd indices
$\ik_E\,=\,\beta_E\,=\,1/p,$
and therefore by the theorem 1.2 $M(E)\subset{L_p}.$

For the proof of inverse embedding we fix $y\in{L_p}$ and consider the linear
operator $T_yx=\,\fg.$ Then
$$
||T_y||_{L_p\to{L_p}}\,=\,||y||_p\,\,\,,\,\,\,||T_y||_{L_{p\infty}\to{L_{p\infty}}}\,\le\,{C_5||y||_p},$$
where $C_5>0$ does not depend from $y$ [ 11 ]. Consequently, in view of
interpolation property of $E$
$$
||T_y||_{E\to{E}}\,\,\le\,\,{C_6\,||y||_p},$$
that is, $y\in{M(E)}$ and $||y||_{M(E)}\,\le\,{C_6\,||y||_p}.$

The theorem is proved.\vspace{\baselineskip}

Particularly, $L_{pq}\,(p\le{q}\le{\infty})$ is interpolation space
between $L_p$ and $L_{p\infty}$ [ 12, p.142 ]. Therefore, we obtain the following

{\bf Corollary 1.12}. If $1<p<{\infty},p\le{q}\le{\infty},$ then
$$
M(\,L_{pq}\,)\,\,=\,\,L_p.$$

{\bf Remark 1.13}. In case $q=\infty$ the last statement was proved in [11].

$$\,\,$$

Let us introduce a quantity connected with the boundedness of tensor product
in the symmetric spaces.
We define for a symmetric space $E,\,\,m\in{\bf{N}}$ and $\ik=(\ik_k)_{k=1}^m\in{\bf{R}^m}$
$$
{\cal K}_{E}^m(\ik)\,=\,\,\sup{||\sum\nolimits_{k=1}^m{\ik_ky_k}||_{E}},$$
where the supremum is taken over all $y\ge{0},\,||y||_E=1$ and all
collections of disjointly supported functions $y_k,$ such that
$$
n_{y_k}(s)\,=\,\frac{1}{m}\,n_y(s),\,s>0.$$

{\bf Theorem 1.14}. Let $E$ be a symmetric space on $I$ with an order subcontinuous norm.
Then $B$ is the bounded operator from $E\times{E}$ into $E(I\times{I})$
if and only if there exists $C>0$ such that for all $m\in{\bf{N}},{\ik=(\ik_i)_{i=1}^m}\in{\bf{R}}^m$
$$
{\cal K}^m_E (\ik)\,\le\,{C||\hi||_E}\eqno{(14)}$$

{\bf Proof}. If $y\in{E},||y||_E=1$ and $x=x(s)$ is such as in (9),then functions
$$
B(x,y)(s,t)\,=\,\hi(s)\,y(t)\,\,(s,t\in{[0,1]})$$
and
$$
v(t)\,=\,\qr(t)\,\,(t\in{[0,1]})$$
( $y_k$ are mutually disjoint and $n_{y_k}(u)=\frac{1}{m}n_y(u),\,u>0$ ) are equimeasurable. 
Consequently, the inequality (14) is equivalent to the following:
$$
||\fg||_E\,\,\le\,\,{C\,||x||_E},\eqno{(15)}$$
where $||y||_E=1$ and the function $x$ belongs to the class (9). Hence 
boundedness of the operator $B$ from $E\times{E}$ into $E(I\times{I})$
 implies (14).

Conversely, let us assume that the inequality (14) ( or (15) ) is satisfied. 
Since the norm of $E$ is order subcontinuous, then the inequality (15) 
can be extended, at first, on countably-value and, next, on arbitrary
functions from $E$ as in proof of the theorem 1.5.
\def \st {\phi_{p,{\alpha}}}
\def \uy {||x||_{p,{\alpha},{q}}}
\def \tu {L_{pq}(\log{L})^{\ik}}
\def \yz {\theta}
\def \zx {L_{p\infty}(\log{L})^{-1/p}}
\def \xw {L_{p\infty}}
\def \wz {I\times{I}}
\def \za {L_{pr}\times{L_{pq}}}
\def \wt {L_{p,2p}(\log{L})^{-1/{2p}}}

\begin{center}
$$\,\,$$

2. {\em The tensor product in the Lorentz spaces} $L_{pq}$
\end{center}

Let $1<p<{\infty},\,
1\le{q}\le{\infty},\,{\ik}\in{\bf{R}},\,\st(u)=\,u^{1/p}\ln^{\ik}(e/u).$
The Lorentz-Zygmund space $\tu$ consists of all measurable functions on $I$
 ( or $\wz$ ) for which
$$||x||_{p,\ik,q}\,=\left\{
\begin {array} {ll}
\left\{\int\limits_0^1(x^*(u)\st(u))^q\,{du}/u\right\}^{1/q}, & \mbox{if } q<{\infty} \\
\sup\limits_{0<u\le{1}}(x^*(u)\st(u)),                            & \mbox{if } q=\infty
\end {array} \right.$$
is finite. As in case $L_{pq}$ $\uy$ is not a norm. However using the 
Minkovski inequality ( see, for example, [16, p. 38] ) it can be verified
that this functional is equivalent to the norm $||x||_{p,{\ik},q}'\,=\,||x^{**}||_{p,{\ik},q}.$
Moreover, the fundamental function of space $\tu$ is equivalent to the 
function $\st(u).$ For $\ik=0$ the Lorentz-Zygmund spaces are the standard Lorentz
spaces $L_{pq}$, while for $p=q<\infty$ they are the Orlicz spaces
$L_p(\log{L})^{{\ik}p}.$ For $p=q=\infty,$ they produce the exponential classes
of Zygmund ( see [17] for further details ). The Lorentz-Zygmund spaces 
were studied of a number of authors ( except [17] see, for example, [18] ). \vspace{\baselineskip}

If $E$ is a symmetric space, then we use $E^0$ to denote the norm closure of $ L_{\infty}$ in $E.$

{\bf Theorem 2.1}. Let $1<p<\infty,\,p\le{r}\le{q}\le{\infty},\,\ik=r^{-1}-p^{-1}.$
The operator $B$ is bounded operator from $\za$ into the space $\left(\tu\right)^0\,(\wz)$ 
 ( if $p<r<q=\infty$ ) and into $\tu\,(\wz)$ ( otherwise ).

We need for the proof of this theorem two lemmas. The first of them is
proved immediately using straightforward calculations.

{\bf Lemma 2.2}. Let $1<p<\infty,\,\ik\in{\bf{R}},\,\psi_{p,\ik}=1/{\st}.$

1) If $C=C(p,\ik)>0$ is large enough, then for $v>C$ the distribution function
of $\psi_{p\ik}$ is equivalent to the function
$
N_{p,\ik}(v)\,=\,v^{-p}\ln^{{-p}\ik}v.$

2) If $\ik_i<1/p\,\,(i=0,1),$ then the functions $\psi_{p,\ik_0}\otimes{\psi_{p,\ik_1}}$ 
and $\psi_{p,\ik_0+\ik_1-1/p}$ are equimeasurable.\newpage

The second lemma can be proved using standard interpolation techniques ( we
shall denote by $[X_0,X_1]_{\yz}$ the complex interpolation spaces, as defined, 
for example, in [13] or [19] ). In detail see [11]. 

{\bf Lemma 2.3}. Let $1<p_1\le{p_0}<\infty,\,-\infty<\ik_1\le{\ik_0}<\infty,\,1\le{q_0}\le{q_1}\le{\infty},\,0\le{\yz}\le{1}$ and
$$
\max\,(\,p_0-p_1\,,\,\ik_0-\ik_1\,,\,q_1-q_0\,)\,>\,0\eqno{(16)}$$

Then the space
$$
\left[\,L_{p_0,q_0}(\log{L})^{\ik_0}\,\,,\,\,L_{p_1,q_1}(\log{L})^{\ik_1}\,\right]_{\yz}$$
is isomorphic to $\left(\,L_{pq}(\log{L})^{\ik}\,\right)^0\,,$ if $q_0=q_1=\infty$
and $0<\yz\le{1}$ or $q_1=\infty$ and $\yz=1$, and to $L_{pq}(\log{L})^{\ik}\,,$ otherwise.
The numbers $p,q,\ik$ are defined by the following way:
$$
\frac1p\,=\,\frac{1-\yz}{p_0}\,+\,\frac{\yz}{p_1}\,,\,\frac1q\,=\,\frac{1-\yz}{q_0}\,+\,\frac{\yz}{q_1}\,,\,\ik\,=\,(1-\yz)\ik_0\,+\,{\yz}{\ik_1}.$$

{\bf Remark 2.4}. The inequality (16) implies the embedding
$$
L_{p_0,q_0}(\log{L})^{\ik_0}\,\subset{\,L_{p_1,q_1}(\log{L})^{\ik_1}}$$
that is used by proof of lemma 2.3. However this embedding is true even when (16) is not satisfied [17].

{\bf Proof of the theorem 2.1}. Let us define numbers
$$
s\,=\,\frac{pq}{r}\,,\,\yz\,=\,1-\frac{p}{r}\eqno{(17)}$$
Then $p\le{s}\le{\infty},0\le{\yz}\le{1},p\le{r}\le{\infty},s\le{q}\le{\infty}$ and
statement of theorem is true for extreme values $r$ and $q$, that is, 
$$
B\,:\,L_{pp}\times{L_{ps}}\,\rightarrow\,{L_{ps}(\wz)}\eqno{(18)}$$
$$
B\,:\,L_{p\infty}\times{L_{p\infty}}\,\rightarrow\,{\zx}\eqno{(19)}$$

In fact, (18) follows from the O'Neil theorem ( or from the corollary 1.12 of
present paper ). In view of the Closed Graph Theorem for  proof (19)
it is enough to show that $B(x,y)$ belongs to $\zx\,(\wz),$ if 
$x\in{\xw},y\in{\xw}.$ Lemma 2.2 and definition of $\zx$ imply that this
is true for the functions $x_0(u)=\,u^{-1/p}\,,\,y_0(v)=\,v^{-1/p}$. 
Since these functions have the greatest rearrangement in $L_{p\infty}$, then (19) is proved.

For $\yz$ from (17) we use now the complex method of interpolation. At first,
we remind its following property ( see [19] or [13, p. 125-126] ).

If $(X_0,Y_0)\,,\,(Y_0,Y_1)$ and $(Z_0,Z_1)$ are arbitrary Banach pairs, $T$ is an 
bilinear operator from $X_i\times{Y_i}$ into $Z_i\,(i=0,1),$ then $T$ is bounded operator from 
$[X_0,X_1]_{\yz}\times{[Y_0,Y_1]_{\yz}}$ into $[Z_0,Z_1]_{\yz}\,(0\le{\yz}\le{1}).$

The statement of theorem 2.1 follows now from last remarks and lemma 2.3.\vspace{\baselineskip}

For proof of following result about the exactness of theorem 2.1. 
it is sufficient to consider operator $B$ on the pairs of 
functions of type
$$
\psi_{p\ik}(u)\,=\,u^{-1/p}ln^{-\ik}(e/u)$$
and to use, next, lemma 2.2.

{\bf Theorem 2.5}. If $1<p<\infty$ and $p<r\le{q}\le{\infty},$then the
operator $B$ does not act from $\za$ into the space $L_{pq}(\log{L})^{\beta}\,(\wz)$
for any $\beta>r^{-1}-p^{-1}.$\vspace{\baselineskip}

{\bf Remark 2.6}. If $r=q=\infty,$ then the following assertion is true:  $\zx$ is
the smallest among symmetric spaces $E$ such that $B$ is the bounded operator
from $L_{p\infty}\times{L_{p\infty}}$ into $E(\wz)$ ( see lemma 2.2 and
proof of theorem 2.1 ).\vspace{\baselineskip} 

We get from O'Neil theorem and theorem 2.2 following

{\bf Corollary 2.7}. Let $1<p\le{r}\le{q}<\infty, r^{-1}+q^{-1}-p^{-1}\,\ge{0}.$
If $\ik=\,r^{-1}-p^{-1},\,s^{-1}=\,r^{-1}+q^{-1}-p^{-1},$ then $B$ is
the bounded operator from $\za$ into the space
$$
E(p,r,q)(\wz)\,=\,\,L_{ps}\bigcap{\tu\,(\wz)}.$$

{\bf Remark 2.8}. The space $E(p,q,r)$ is not equal neither $L_{ps}$ nor $\tu,$ if $p<r.$
Assume, for example, that $r=q=2p.$ Then $B$ acts from $L_{p,2p}\times{L_{p,2p}}$
into $E(p,2p,2p)(\wz)=\,\xw\bigcap{\wt\,(\wz)}.$  

On the one hand, it is easy to sheck that the function $\psi_{p,0}(u)=\,u^{-1/p}$
belongs to $\xw$ and it does not belong to $\wt$.

On the other hand, if we assume that, vice versa,
$$
\wt\subset{\xw},$$
then from the inequality for fundamental functions of these spaces 
we obtain false inequality:
$\ln(e/t)\,\le\,{C}$
( $C>0$ does not depend on $t\in{(0,1]}$ ).\vspace{\baselineskip}

{\bf Remark 2.9}. As we remarked in beginning of this paper the case
$$
\frac{1}{r}+\frac{1}{q}-\frac{1}{p}\,<\,0$$
is not considered in the O'Neil theorem. The theorems 2.1 and 2.4 complete
it, characterising the image $B(L_{pr},L_{pq})$ in this case to within
any positive power of logarithm.

{\bf References}
\begin {enumerate}

\item M.Milman, Some new function spaces and their tensor products, Notes
de Mat., 20 (1978) 1-128.
\item J.Lindenstrauss and L.Tzafriri, Classical Banach spaces 2, Function
spaces, Springer, Berlin 1979.
\item W.B.Johnson, B.Maurey, G.Schechtman and L.Tzafriri, Symmetric 
structures in Banach spaces, Mem. AMS 217 (1979) 1-298.
\item N.L.Carothers, Rearrangement invariant subspaces of Lorentz function
spaces, Isr. J. Math. 40, No.3-4 (1981) 217-228.
\item C.Corduneanu, Integral equations and stability of Feedback systems,
 Academic Press, New York 1973.
\item M.Milman, Tensor products of function spaces, Amer. Math. Soc. 82,
No.4 (1976) 626-628.
\item M.Milman, Embeddings of Lorentz-Marcinkiewicz spaces with mixed
norms, Anal. math. 4, No.3 (1978) 215-223.
\item M.Milman, Embeddings of $L(p,q)$ spaces and Orlicz spaces with
mixed norms, Notes de Mat. 13 (1977) 1-7.
\item R.O'Neil, Integral transforms and tensor products on Orlicz spaces
and $L(p,q)$ spaces, J. d'Analyse Math. 21 (1968) 1-276.
\item R.O'Neil, Convolution operators and $L(p,q)$ spaces, Duke Math. J.
30 (1963) 129-142.
\item S.V.Astashkin, On bilinear multiplicative operator [in Russian],
In: Theory of functions of several variables, Jaroslavl' (1982) 3-15.
\item S.G.Krein, Yu.I.Petunin and E.M.Semenov, Interpolation of Linear
Operators [in Russian], Nauka, Moskow 1978.
\item J.Bergh and J.Lofstrom, Interpolation Spaces. An Introduction
[Russian translation], Mir, Moscow 1980.
\item H.P.Rosenthal, On a theorem of J.L.Krivine concerning block 
finite-representablity of $l^p$ in general Banach spaces, J. Funct.
Anal. 28 (1978) 197-225.
\item L.V.Kantorovicz and G.P.Akilov, Functional Analysis [in Russian],
Nauka, Moscow 1977.
\item A.Zygmund, Trigonometric Series.V.1 [Russian translation], Mir,
Moskow 1965.
\item R.Sharpley, Counterexamples for classical operators on 
Lorentz-Zygmund spaces, Studia Math. 68, No.2 (1980) 141-158.  
\item B. Jawerth and M.Milman, Extrapolation theory with application,
Mem. AMS 89,No.440 (1991) 1-82.
\item A.P.Calderon, Intermediate spaces and interpolation, the complex
method, Studia Math. 24 (1964) 113-190.
\end{enumerate}
\end{document}